\documentclass[11pt]{amsart}
\usepackage{amscd,amssymb,palatino}
\usepackage[utf8]{inputenc}
\usepackage[english]{babel}
\usepackage{caption}
\usepackage{etex}
\usepackage[leqno]{amsmath}
\usepackage{amssymb}
\usepackage{color}
\usepackage{shadow}
\usepackage{epsfig}
\usepackage{epic}
\usepackage{eepic}
\usepackage{graphics}
\usepackage{graphicx}
\usepackage{psfrag}
\usepackage{calc}
\usepackage{tikz}
\usepackage[dvipsnames,prologue,table]{pstricks}
\usepackage{pstcol}
\usetikzlibrary{arrows,decorations,patterns,positioning,automata,shadows,fit,shapes,calc,decorations.markings,backgrounds,scopes,decorations.text}

\usepackage{tipa}

\usepackage{fancyhdr}
\fancypagestyle{plain}{\fancyhf{}
\fancyfoot[C]{{\large\sf\bfseries\thepage}}}
\usepackage{calc}
\newcommand\1{\lower 9pt\hbox{\underbar{}}}
\numberwithin{equation}{section}

\newtheorem {Theorem}                   {Theorem}

\newtheorem {Lemma}[equation]           {Lemma}

\theoremstyle{definition}
\newtheorem {Definition}[equation]{Definition}
\newtheorem {Remark}[equation]          {Remark}

\newcommand{\pr} {\smallskip\noindent{\bf Proof\,\,}}
\newenvironment{Proof}  {\pr}{\hspace*{\fill}\qed\\}

\begin{document}

\title{Convexity of the moment map image for torus actions on $b^m$-symplectic manifolds}

\author{Victor W. Guillemin}
\author{Eva Miranda}
\thanks{{ E. Miranda  is supported by the Catalan Institution for Research and Advanced Studies via an ICREA Academia Prize 2016, a Chaire d'Excellence de la Fondation Sciences Math\'{e}matiques de Paris and partially supported  by the grants reference number MTM2015-69135-P (MINECO/FEDER) and reference number 2014SGR634 (AGAUR).This work is supported by a public grant overseen by the French National Research Agency (ANR) as part of the \emph{\lq\lq Investissements d'Avenir"} program (reference: ANR-10-LABX-0098).}}
\author{Jonathan Weitsman}
\thanks{J. Weitsman was supported in part by NSF grant DMS 12/11819}
\address{Department of Mathematics, MIT, Cambridge, MA 02139}
\email {vwg@math.mit.edu}
\address{{Department of Mathematics}, Universitat Polit\`{e}cnica de Catalunya and BGSMath, Barcelona, Spain \\ CEREMADE at Universit\'{e} de Paris Dauphine,  IMCCE at Observatoire de Paris and IMJ at Université de Paris Diderot, Postal address: 77, avenue Denfert Rochereau
75014, Paris, France}
\email{eva.miranda@upc.edu, Eva.Miranda@obspm.fr}
\address{Department of Mathematics, Northeastern University, Boston, MA 02115}
\email{j.weitsman@neu.edu}
\thanks{\today}

\begin{abstract} We prove a convexity theorem for the image of the moment map of a Hamiltonian torus action on a $b^m$-symplectic manifold.
\end{abstract}
\maketitle

\section{Introduction}

The purpose of this paper is to prove a convexity theorem for the image of the moment map of a Hamiltonian torus action on a $b^m$-symplectic manifold.
$b^m$-symplectic manifolds are Poisson manifolds where the Poisson structure is invertible on the complement of a hypersurface $Z$, and has a singularity of order $m$ on $Z,$ where $m \geq 1$ is an integer. They are a generalization of $b$-symplectic (or, log-symplectic) manifolds and were introduced in the thesis of G. Scott \cite{scott}.

A convexity theorem for the moment image of a Hamiltonian torus action on a $b$-symplectic manifold was proved in \cite{gmps2}.  In this case, the moment image is governed by the singularity of the moment map in the neighborhood of $Z$, encoded in the {\em modular weight} (see \cite{gmps}).  If this modular weight is nonzero, the image is not only convex, but on each component of the complement of $Z,$ has the form of a product of a convex polytope with a ray or the real line, possibly modified by symplectic cutting.

We show that the moment image has a similar form where $m > 1.$  The argument requires some care, since in this case the moment map has an asymptotic series near $Z$ which involves $m$ modular weights.  The leading term in this series gives the moment map a form as in \cite{gmps2}.  We show that the subleading terms preserve this structure.

As a technical device to state our convexity theorem, we use the desingularization of {$b^m$-}symplectic forms in \cite{gmw1}.  This simplifies the statements since we can deal with families of compact manifolds, with compact moment images.

In a companion paper, we use the convexity theorem to study formal geometric quantization of $b^m$-symplectic manifolds equipped with Hamiltonian torus actions.  The form of the moment image will give rise to a very simple asymptotics of this quantization.

{\textbf{Acknowledgements:} We thank C\'{e}dric Oms for carefully reading a first version of this article.}
\section{$b^m$-Manifolds}

Let $M$ be a compact manifold, and let $f\in\mathcal{C}^\infty(M)$ have a transverse zero at  a hypersurface $Z\subset M$. Let $m$ be a positive integer, the $m$-germ of $f$ at $Z$ gives rise to a sheaf, and therefore a vector bundle, whose sections are given by,

$$\Gamma(^{b^m}TM)=\{v \in \Gamma(TM): v \quad \text{vanishes to order $m$ at $Z$} \}.$$

\noindent  By considering sections of the wedge powers {${\Lambda^k} (^{b^m}{T^*M})$} we obtain a complex $(^{b^m}\Omega^k(M), d)$ of differential forms with singularities at $Z$.

The  cohomology associated to this complex is given by the following:

\begin{Theorem}[\textbf{$b^m$-Mazzeo-Melrose}, \cite{scott}]\label{thm:Mazzeo-Melrose}
\begin{equation}\label{eqn:Mazzeo-Melrosebm}
  ^{b^m}H^p(M) \cong H^p(M)\oplus(H^{p-1}(Z))^m.
\end{equation}
\end{Theorem}

This theorem comes from a Laurent expansion in a neighborhood of $Z$ which for $2$-forms has the form:

\begin{equation}\label{eqn:newlaurent}
{\omega = \sum_{j = 1}^{m}\frac{df}{f^j} \wedge \pi^*(\alpha_{j}) + \beta}
\end{equation}
\noindent where $\alpha_j$ are closed one forms on $Z$, $\beta$  is a  closed 2-form on $M$, and $\pi:U\longrightarrow Z$ is the projection.

A two  $b$-form $\omega\in {^{b^{m}}\Omega}^2(M)$ is $b^m$-symplectic if it is closed and nondegenerate as an element of $\Lambda^{2}(  {^{b^{m}}T^*M})$.

\begin{Theorem}\label{newmoser}
  There exists a neighborhood $U=(-\epsilon, \epsilon)\times Z$ of $Z$ and a diffeomorphism $\psi: U\longrightarrow U$ preserving $Z$ and the $m$-germ of $f$ such that

\begin{equation}\label{eqn:newlaurent}
{\psi^*(\omega) = \sum_{j = 1}^{m}\frac{df}{f^j} \wedge \pi^*(\alpha_{j}) + \pi^*(\beta)}
\end{equation}
 where
$\pi:U\longrightarrow Z$ is the projection and $\alpha_i$, and $\beta$ are closed.
\end{Theorem}

\begin{Proof}

  We know by \cite{guimipi} proposition 3.1 that

  {$$\omega = \sum_{j = 1}^{m}\frac{df}{f^j} \wedge \pi^*(\alpha_{j}) + \beta.$$}

  Let $i:Z\longrightarrow U$ be the restriction. Then

  $$\pi^* i^*\alpha_i-\alpha_i= da_i$$
  $$\pi^* i^*\beta-\beta= db$$

  and for $\epsilon$ sufficiently small, letting

  {$$\omega_0= \sum_{j = 1}^{m}\frac{df}{f^j} \wedge \pi^* i^*(\alpha_{j}) + \pi^*i^*(\beta)$$}

  \noindent that $\omega_t= t\omega_0+ (1-t) \omega$ is symplectic for all $t\in [0,1]$.

  Thus

 $$ \frac{d}{dt}\omega_t= d b_t$$

 \noindent where $$b_t=\iota_{v_t}\omega_t$$

 \noindent for $v_t$ a section of $^{b^m} TM$.

 By shrinking $\epsilon$ we may chose $v_t$ compactly supported on $U$  and $\vert\vert v_t \vert\vert <\delta$ for  any fixed $\delta >0$. Thus the existence theorem for ODE shows that $v_t$ integrates to a family of diffeomorphisms, $\psi_t$, vanishing to order $m$ on $Z$ such that $\psi_t^*(\omega_0)=\omega$.

\end{Proof}

\section{Torus actions on $b^m$-manifolds} Now let us assume a torus $T$ acts  on $(M,Z,f)$  preserving $\omega$.
We denote $^{b}\mathcal{C}^\infty(M)=\{ a\log\vert f\vert+ g, g\in \mathcal{C}^\infty(M)\}$  the space of smooth functions with logarithmic singularity at $Z$ and write $^{b^m} \mathcal{C}^\infty(M)=\mathcal{C}^\infty(M)\oplus \left(\oplus_{j=1}^{m-1} f^{-j}\mathcal{C}^\infty(M)\right )\oplus ^{b} \mathcal{C}^\infty(M)$.

\begin{Definition}

The action of $T$ is \textbf{Hamiltonian} if there exists a moment map

$\mu\in^{b^m} \mathcal{C}^\infty(M)\otimes \mathfrak{t}^*$ with $$\langle d\mu, X \rangle= \iota_{X^M}\omega$$

\noindent for any $X\in\mathfrak{t}$ where $X^M$ is the fundamental vector field generated by $X$ on $M$.

\end{Definition}

We will prove the following

\begin{Lemma}
There exists a neighborhood $U=Z\times (-\epsilon, \epsilon)$ where the moment map $\mu: M\longrightarrow\mathfrak{t}^*$ is given by
{$$\mu= a_1\log {\vert f\vert}+\sum_{i=2}^{m}a_i \frac{f^{-(i-1)}}{i-1}  + \mu_0$$}
\noindent  with $a_i\in\mathfrak{t^0}_L$, and $\mu_0$ is the moment map for the $T_L$-action on the symplectic leaves of the foliation.

\end{Lemma}

\begin{Proof}

Note that theorem \ref{newmoser} holds equivariantly so we may assume that $\omega$ can be written \begin{equation}\label{eqn:newlaurent2}
{\omega = \sum_{j = 1}^{m}\frac{df}{f^j} \wedge \pi^*(\alpha_{j}) + \pi^*(\beta)}
\end{equation}
 where
$\pi:U\longrightarrow Z$ is the projection and $\alpha_i$, and $\beta$ are closed and $T$-invariant.

%

The moment map $\mu$ therefore has the form,\begin{equation}
\label{eqn:momentmapeqn}
{\mu= a_1\log {\vert f\vert}+\sum_{i=2}^{m}a_i \frac{f^{-(i-1)}}{i-1}  + \mu_0}
\end{equation}
\noindent where  $\left< a_i, X \right>=\alpha_i(X^M)$ and  $\mu_0$ is the moment map for the action of $T$ on the regular Poisson manifold $Z$.

\end{Proof}
The form $\alpha_m$ is nowhere vanishing and determines the symplectic foliation of $Z$. We now make the following two assumptions,

\begin{itemize}
  \item \textbf{Assumption 1}
There exists $\xi \in \mathfrak{t}$ such that $a_m(\xi)\neq 0$.  Remark: Without loss of generality, we may assume also that $\xi$ generates a circle  subgroup $S^1_{\xi}\in T$.
  \item \textbf{Assumption 2} The foliation given by $\alpha_m$ has a compact leaf $L$ . Thus $(L,\beta)$ is a compact symplectic manifold acted on by the Hamiltonian torus action of $T_L=T/S^1_{\xi}$, with moment map ${\mu_0}_{\vert\mathfrak{t}_L }$. This would follow from integrality of $\omega$; see \cite{gmw4}.

\end{itemize}

\section{Desingularization}

Recall from \cite{gmw1},

\begin{Theorem}\label{thm:deblogging}
Given a $b^m$-symplectic structure $\omega$ on a compact manifold $M^{2n}$  let $Z$ be its critical hypersurface.
\begin{itemize}
\item If $m$ is even, there exists  a family of symplectic forms ${\omega_{\epsilon}}$ which coincide with  the $b^{m}$-symplectic form
    $\omega$ outside an $\epsilon$-neighborhood of $Z$ and for which  the family of bivector fields $(\omega_{\epsilon})^{-1}$ converges in
    the $C^{m-1}$-topology to the Poisson structure $\omega^{-1}$ as $\epsilon\to 0$ .
\item If $m$ is odd, there exists  a family of folded symplectic forms ${\omega_{\epsilon}}$ which coincide with  the $b^{m}$-symplectic form
    $\omega$ outside an $\epsilon$-neighborhood of $Z$.

\end{itemize}

\end{Theorem}

If a torus acts, this family of forms $\omega_{\epsilon}$ may be chosen equivariantly.

\begin{Remark}
Observe that even if the initial action is Hamiltonian in the $b^m$-sense, the desingularized action need not be Hamiltonian in the standard sense because its moment map might be circle valued.

\end{Remark}

\section{The local convexity theorem}

\begin{Theorem}
  Let $Z_i$ be a connected component of $Z$. Then there exists a neighborhood $(-\epsilon, \epsilon)\times Z_i$ where the image of the moment map for the $T$-action on the desingularized family $(M,\omega_{\epsilon})$ is

  \begin{itemize}
    \item  $\Delta_i\times (-a_{\epsilon},a_{\epsilon})$ for even $m$.
    \item $\Delta_i\times (-a_{\epsilon},a_{\epsilon})/\psi$ for odd $m$.

  \end{itemize}

  \noindent where $\Delta_i$ is the image of the moment map for the $T_{L_i}$-action on $L_i$,  a symplectic leaf on $Z$ and

  \begin{itemize}
    \item $a_{\epsilon}\to\infty$ as $\epsilon\to 0$.
    \item $\psi: (-a_{\epsilon},a_{\epsilon})\mapsto (-a_{\epsilon},a_{\epsilon})$ is the involution $x\mapsto -x$.
  \end{itemize}
\end{Theorem}

\begin{Remark} This implies that the image of the moment map is locally convex.
\end{Remark}
\begin{Proof} Recall the moment map is given by the expression (\ref{eqn:momentmapeqn}) with $a_i\in \mathfrak{t}^*$ constant.We claim that in (\ref{eqn:momentmapeqn})  $\left < a_i, \xi \right >=0$ for all $\xi\in \mathfrak{t}_L$.
To see this, suppose there exists $\xi\in \mathfrak{t}_L$ such that $\left < a_i, \xi \right > \neq 0$. Since $T_L$ is a torus action on a compact Hamiltonian $T_L$-space, $L$, it must have a fixed point $p\in Z$. At this fixed point $\xi^M_{\vert p}=0$ so
$\left < a_i, \xi \right >=\alpha_i(\xi^M_{\vert p})=0$. Thus $\left < a_i, \xi \right >=0$.
\end{Proof}

\section{Global convexity theorem}

In this section we prove,

\begin{Theorem}
  Let $(M, Z, f)$ be a $b^m$-symplectic manifold. Let $M\setminus Z= \bigsqcup_{i=1}^r M_i$. Then the image of the moment map for the desingularized symplectic form  $\omega_{\epsilon}$ on $\bar{M_i}$ is given by either:

    \begin{enumerate}
    \item  A product $\Delta\times [-a_{\epsilon},a_{\epsilon}]$, where $a_\epsilon \to \infty$ as $\epsilon \to 0;$
    or,
    \item A convex polytope which has the form of a product $\Delta\times [0,a_{\epsilon}]$ in the neighborhood of $Z$; in other words, a polytope of the form

    $\Delta\times [-a_{\epsilon},a_{\epsilon}]\cap H_1\cap \dots \cap H_n$.

  \end{enumerate}

  \noindent Here $\Delta$ is the image of the moment map for the $T_L$-action on $L$ and $H_1,\dots, H_n$ are half-spaces. In particular the image polytopes $\Delta_i$ coincide.

  \end{Theorem}

  \begin{Proof}By the local convexity theorem, we know that if $\partial M_i = Z_i\bigsqcup Z_{i+1}$, then in a neighborhood of $Z_i$, respectively $Z_{i+1}$ the image is of the form $\Delta_i\times [-a_{\epsilon},c]$ resp $\Delta_{i+1}\times [c',a_{\epsilon}]$ where $a_{\epsilon}\to \infty$ as $\epsilon \to 0$ and we may take $c$ and $c'$ positive.

 Let $\nu=\mu_{\vert M_i}$ denote the restriction of $\mu$ to $M_i$ and let $P=\nu(M_i)$ be the image of the moment map. By the convexity theorem for torus actions on symplectic manifolds, the image $P$ of the moment map restricted to $M_i$ is convex.

  Therefore $P\cap \nu^{-1}( [-a_{\epsilon}, -c])=[-a_{\epsilon}, -c] \times \Delta_i$ and $P\cap \nu^{-1}([c', a_{\epsilon}])=[c', a_{\epsilon}]\times \Delta_{i+1}$ for all $c, c'$ sufficiently close to $a_{\epsilon}$.

Let us see that this implies $P\cap \nu^{-1}([-a_{\epsilon}, a_{\epsilon}])=[-a_{\epsilon}, a_{\epsilon}] \times \Delta_i.$

We first prove $P\cap [-a_{\epsilon}, a_{\epsilon}]\subset [-a_{\epsilon}, a_{\epsilon}] \times \Delta_i$. Suppose $p\in P \cap \nu^{-1}([-a_{\epsilon},a_{\epsilon}])$ but $p\notin [-a_{\epsilon},a_{\epsilon}]\times \Delta_i$; then $p=(t,\tau)$ where $t\in [-a_{\epsilon}, a_{\epsilon}]$ with $\tau \notin \Delta_i$.

Let $\sigma$ be a point of $\Delta_i$ at minimal distance from $\tau$. The line $\ell$ connecting $(t, \tau)$ and $(-\alpha, \sigma)$ must lie in $P$. Thus for a sufficiently close to $-\alpha$, $\ell\cap[-\alpha, a]\times \mathbb R^{d-1}$ must lie in $[-\alpha, a]\times \Delta_i$. Hence $\ell$ must intersect the plane $(a, x)_{x\in \mathbb R^{d-1}}$ at a point $\sigma'\in \Delta_i$ closer to $\tau$ than $\sigma$ (see picture below). This is a contradiction.

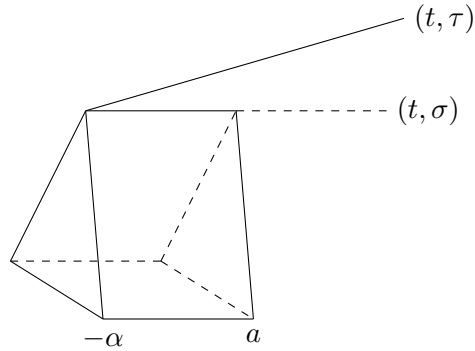
\begin{center}
\begin{tikzpicture}[scale=2]

\draw [-] (0,0,0) -- (1,0,0) node[left] {};
\draw [-] (0,0,0) -- (-1,0,-1) node[left] {};
\draw [-,dashed] (0,0,-1) -- (1,0,0) node[left] {};
\draw [-,dashed] (0,0,-1) -- (-1,0,-1) node[left] {};

\draw [-] (0,0,0) -- (-1/2,1,-1) node[left] {};
\draw [-] (-1,0,-1) -- (-1/2,1,-1) node[left] {};
\draw [-] (1,0,0) -- (-1/2+1,1,-1) node[left] {};
\draw [-,dashed] (0,0,-1) -- (-1/2+1,1,-1) node[left] {};

\draw [-] (-1/2+1,1,-1) -- (-1/2,1,-1) node[left] {};

\draw (0,0,0) node [below] {$-\alpha$};
\draw (1,0,0) node [below] {$a$};

\draw [-] (-1/2,1,-1) -- (2,2,0) node[right]{$(t,\tau)$};
\draw [-,dashed] (-1/2+1,1,-1) -- (-1/2+2,1,-1) node[right]{$(t,\sigma)$};

\end{tikzpicture}
\end{center}
Thus we have proved that  $P\subset [-a_{\epsilon}, a_{\epsilon}] \times \Delta_i$; also $P\subset [-a_{\epsilon}, a_{\epsilon}] \times \Delta_{i+1}$.
Therefore, $P\subset [-a_{\epsilon}, a_{\epsilon}] \times (\Delta_i\cap \Delta_{i+1})$. In particular for $c$ sufficiently close to $a_{\epsilon}$,
$[-a_{\epsilon}, c]\times \Delta_i \subset [-a_{\epsilon}, c]\times (\Delta_i\cap \Delta_{i+1})$ and $[c, a_{\epsilon}]\times \Delta_i \subset [c, a_{\epsilon}]\times (\Delta_i\cap \Delta_{i+1})$ so $\Delta_i \subset (\Delta_i\cap \Delta_{i+1})$ and $\Delta_{i+1} \subset (\Delta_i\cap \Delta_{i+1})$
; thus $\Delta_i=\Delta_{i+1}$. Since $P$ contains all lines $[-a_{\epsilon}, a_{\epsilon}]\times \Delta_i$, it must be equal to $[-a_{\epsilon}, a_{\epsilon}]\times \Delta_i$.
On the other hand, if $\partial M_i$ is connected the image of the moment map must be of the form (1) or (2) and can be written as $\Delta\times [-a_{\epsilon},a_{\epsilon}]\cap H_1\cap \dots \cap H_n$ for some half-spaces $H_1, \dots H_n$.
  \end{Proof}

\end{document}